\documentclass[twocolumn,showpacs,aps,floatfix]{revtex4}

\usepackage{bm}
\usepackage{amsmath}
\usepackage{amssymb}
\usepackage{latexsym} 
\usepackage{amsfonts} 
\usepackage{epsfig}
\usepackage{psfrag}

\newcommand{\nn}{\nonumber\\} 
 
\newcommand{\f}[1]{\mbox{\boldmath$#1$}}

\newcommand{\bea}{\begin{eqnarray}}
\newcommand{\ea}{\end{eqnarray}}
\newcommand{\eea}{\end{eqnarray}}

\newcommand{\commute}[2]{{\left[#1\,,#2\right]}}
\newcommand{\norm}[1]{\left\| #1 \right\|}
\newcommand{\abs}[1]{\left|#1\right|}

\begin{document}
  
\title{A fixed point iteration for computing the matrix logarithm} 

\author{Gernot Schaller$^1$}

\affiliation{$^1$Institut f\"ur Theoretische Physik, 
Technische Universit\"at Dresden, 01062 Dresden, Germany}  

\begin{abstract}
In various areas of applied numerics, the problem of calculating the logarithm of a matrix $A$ emerges.
Since series expansions of the logarithm usually do not converge well for matrices far away from the
identity, the standard numerical method calculates successive square roots.
In this article, a new algorithm is presented that relies on the computation of successive matrix 
exponentials.
Convergence of the method is demonstrated for a large class of initial matrices and 
favorable choices of the initial matrix are discussed.
\end{abstract}


\maketitle
%%%%%%%%%%%%%%%%%%%%%%%%%%%%%%%%%%%%%%%%%%%%%%%%%%%%%%%%%%%%%%%%%%%%%%%%%%%%%%%
%%%%%%%%%%%%%%%%%%%%%%%%%%%%%%%%%%%%%%%%%%%%%%%%%%%%%%%%%%%%%%%%%%%%%%%%%%%%%%%
%%%%%%%%%%%%%%%%%%%%%%%%%%%%%%%%%%%%%%%%%%%%%%%%%%%%%%%%%%%%%%%%%%%%%%%%%%%%%%%
%%%%%%%%%%%%%%%%%%%%%%%%%%%%%%%%%%%%%%%%%%%%%%%%%%%%%%%%%%%%%%%%%%%%%%%%%%%%%%%

\section{Introduction}\label{Sintro}

Calculating the logarithm of a quadratic matrix $A$ may be a difficult problem.
Of course, for problems of moderate size it is usually not prohibitive to calculate 
its logarithm (if existent) by complete diagonalization \cite{press1994}.
However, there may exist several arguments against this method:
In some cases, $A$ might not be diagonalizable. 
Also, if $A(t)$ and $\ln A(t)$ are known for some given time $t$, one might want
to have an efficient solution for $\ln A(t+dt)$ where $dt$ is small.
The approximate calculation of eigenvalues and eigenvectors of $A$ with Arnoldi 
methods \cite{arpack1998} would enable one to follow a time-dependent spectrum
efficiently.
Unfortunately, these methods only yield a part of the spectrum.
In such cases, direct diagonalization is probably not the best method to
compute $\ln A$.

If the matrix $A$ is near the identity matrix, 
one may truncate the Taylor series expansion of the logarithm
\bea
\ln[\f{1}-(\f{1}-A)] = -\sum_{n=1}^\infty \frac{(\f{1}-A)^n}{n}\,.
\eea
However, if this is not the case, convergence may become extremely slow or even fail, 
such that such a series expansion is of little practical use.
The standard resolution to this problem is to bring $A$ near the identity by repeatedly
computing its square root
\bea
\ln A = 2^k \ln A^{1/2^k}\,.
\eea
If the square-roots are only approximated, this can be adapted to
an efficient method \cite{cheng2001a,cardoso2003a}.

The present work undertakes a different step to decrease the computational 
burden.
Instead of calculating successive square roots, a fixed-point iteration 
is presented that requires the calculation of successive exponentials.
The article is organized as follows: After discussing the fixed point iteration scheme
in \ref{Sfp}, in section \ref{Soptimization} further improvements are discussed.
In section \ref{Salgorithm}, a numerical implementation is discussed, and its performance
is analyzed in \ref{Sperformance}.

%%%%%%%%%%%%%%%%%%%%%%%%%%%%%%%%%%%%%%%%%%%%%%%%%%%%%%%%%%%%%%%%%%%%%%%%%%%%%%%
%%%%%%%%%%%%%%%%%%%%%%%%%%%%%%%%%%%%%%%%%%%%%%%%%%%%%%%%%%%%%%%%%%%%%%%%%%%%%%%
%%%%%%%%%%%%%%%%%%%%%%%%%%%%%%%%%%%%%%%%%%%%%%%%%%%%%%%%%%%%%%%%%%%%%%%%%%%%%%%
%%%%%%%%%%%%%%%%%%%%%%%%%%%%%%%%%%%%%%%%%%%%%%%%%%%%%%%%%%%%%%%%%%%%%%%%%%%%%%%

\section{Fixed point iteration scheme}\label{Sfp}

Consider the iteration formula
\bea\label{Eiteration}
X_{n+1} = g(X_n) = A e^{-X_n} - \f{1} + X_n\,.
\eea
If one regards $X_i$ as real numbers, it is quite straightforward to see that 
for any $X_0 > 0$ the above iteration formular converges to
\bea
\lim_{n\to\infty} X_n = \ln A\,.
\eea
For example, it is immediately evident that $X_\infty=\ln A$ is the only fixed point of 
$g(x)$ in \ref{Eiteration}.
In addition, by investigating that $|g'(X_\infty)|=0<1$ one can conclude that this 
fixed point is stable and thus that one has a contractive map which converges to
$\ln \alpha$ for all positive numbers.

Evidently, one can also consider the iteration (\ref{Eiteration})
for matrices $X_n$ and $A$ (assumed to have a well-defined logarithm).
Clearly, one still has fixed points at the logarithms 
$g(\ln(A) + 2 k \pi i) = \ln(A) + 2k\pi i$.
The question is under which conditions these fixed points are attractive, i.e., for which 
matrices $X_i$ the difference (with respect to some norm) to $\ln A$ becomes smaller 
with each iteration.

For simplicity I will restrict myself to an initial matrix $X_0$ that commutes with
$A$.
It is straightforward to show that
\bea
\commute{X_i}{A}=0 \qquad\Longrightarrow\qquad \commute{X_{i+1}}{A}=0\,,
\eea
and therefore if $\commute{X_0}{A}=0$ the iteration (\ref{Eiteration}) defines a
series of mutually commuting matrices.

With inserting $X_i = \ln(A) + \Delta_i$ one has (using $\commute{A}{\Delta_i}=0$)
\bea\label{Edevred}
\Delta_{i+1} = e^{-\Delta_i} + \Delta_i - \f{1}\,.
\eea
From now on I will assume that $\Delta_i$ is a normal matrix, i.e., 
$\commute{\Delta_i}{\Delta_i^\dagger}=0$.
The spectral theorem implies the existence of an orthonormal basis, within which 
$\Delta_i$ has diagonal form (with possibly complex eigenvalues $\lambda_j^{(i)}$).
Therefore, the iteration (\ref{Eiteration}) transforms the eigenvalues of the 
deviation matrix $\Delta_i$ according to (\ref{Edevred}).

If the eigenvalues of $\Delta_i$ are real, one can deduce from 
\bea
e^{-\lambda} + \lambda - 1 \ge 0 \qquad : \qquad -\infty < \lambda < \infty
\eea
that the next deviation matrix $\Delta_{i+1}$ will be positive semidefinite.
Also, one obtains for the operator norm (for a positive semidefinite matrix this is simply 
its largest eigenvalue)
\bea
\norm{\Delta_{i+1}} = \max_j \left[e^{-\lambda_j^{(i)}} + \lambda_j^{(i)} - 1\right]\,.
\eea
Then, one can deduce from
\bea\label{Eexpineq}
e^{-\lambda} + \lambda - 1 &\le& \abs{\lambda} \qquad : \qquad \lambda \ge -1.256\ldots
\eea
that the norm of $\Delta_{i+1}$ will be smaller than the norm of $\Delta_i$ if $\Delta_i$ is 
positive semidefinite.
In other words, for any matrix $X_0 = \ln A + \Delta_0$ with a self-adjoint initial
deviation matrix $\Delta_0$ and $\commute{\Delta_0}{A}=0$ the iteration (\ref{Eiteration}) is 
contractive and will converge to $\ln A$.
Note that (under the precondition that $\commute{A}{\Delta_0}=0$) 
$\Delta_0=\Delta_0^\dagger$ does not imply that $A=A^\dagger$, but
$A=A^\dagger$ implies $\Delta_0=\Delta_0^\dagger$.

\begin{figure}[ht]
\includegraphics[height=7.0cm]{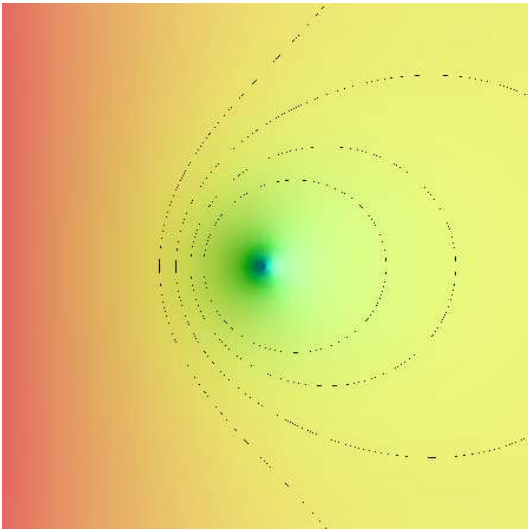}
\caption{\label{Flogarithm}
Density plot of $f(\lambda) = \ln\left(\abs{e^{-\lambda}+\lambda-1}^2\right)-\ln\left(\abs{\lambda}^2\right)$.
The origin is in the center and both real part of $\lambda$ (horizontal axis) and imaginary part 
(vertical axis) range from $-\pi\ldots+\pi$.
The isolines represent values of $f(\lambda)=\{0, 0.5, 1.0, 1.5\}$, respectively, such that the leftmost isoline
borders the region of sure convergence (right, green to blue colours).
In any case, convergence is ensured if the eigenvalues of $\Delta$ are not too far from the 
real axis and if their real part is positive.
}
\end{figure}

Complex eigenvalues $\lambda=x+iy$ (in case of $\Delta_i$ being normal) 
are also transformed according to (\ref{Edevred}).
Demanding that the modulus of all eigenvalues of $\Delta$ should become smaller with each 
iteration, one obtains a region of convergence $V$.
If all eigenvalues of the $\Delta$-matrix are contained within
\bea
V = \{\lambda=x+i y \in C \;:\;\abs{e^{-\lambda}+\lambda-1}^2 \le \abs{\lambda}^2\}
\eea
convergence is assured, see also figure \ref{Flogarithm}.
Note that the ambiguity of the logarithm function does not pose a major problem here, since 
the $\Delta$-matrix can be chosen to represent the difference to any specific branch of the 
logarithm without any difference.

With the Banach fixed point theorem one can then show that the iteration (\ref{Eiteration}) will
converge towards $\ln A$.
In the following, some suitable choices for an initial matrix to the algorithm will be discussed.

%%%%%%%%%%%%%%%%%%%%%%%%%%%%%%%%%%%%%%%%%%%%%%%%%%%%%%%%%%%%%%%%%%%%%%%%%%%%%%%
%%%%%%%%%%%%%%%%%%%%%%%%%%%%%%%%%%%%%%%%%%%%%%%%%%%%%%%%%%%%%%%%%%%%%%%%%%%%%%%
%%%%%%%%%%%%%%%%%%%%%%%%%%%%%%%%%%%%%%%%%%%%%%%%%%%%%%%%%%%%%%%%%%%%%%%%%%%%%%%
%%%%%%%%%%%%%%%%%%%%%%%%%%%%%%%%%%%%%%%%%%%%%%%%%%%%%%%%%%%%%%%%%%%%%%%%%%%%%%%

\section{Algorithmic Optimizations}\label{Soptimization}

Beyond performing scaling transformations $\ln A = \ln(A/\sigma) + \ln(\sigma)\f{1}$, one
has further options to improve the algorithmic performance.

It is evident from (\ref{Eexpineq}) and figure \ref{Flogarithm} that a good initial 
guess for the matrix logarithm may save a lot of computation time.
Such a guess can be made if some bounds on the eigenvalues of $A$ (and hence also
to those of $\ln A$) are known.
Such bounds can for example be cheaply extracted from Gershgorins circle 
theorem (which is especially useful if the matrix $A$ is diagonally dominant) 
or they may be already known from the definition of the problem.
Then an initial matrix with eigenvalues close to $\ln A$ can be constructed from
the linear Taylor approximation that could be optimized for the regime 
$\left[\ln(\lambda_{\rm min}),\ln(\lambda_{\rm max})\right]$.
Assuming an approximately uniform distribution of eigenvalues, one could for example
consider 
\bea
X_0 = \left[\ln\left(\frac{\lambda_{\rm min} + \lambda_{\rm max}}{2}\right)-1\right]\f{1} 
+ \frac{2}{\lambda_{\rm min} + \lambda_{\rm max}} A\,.
\eea
Other choices could include some adapted polynomials of $A$.

In addition to a good guess for an initial matrix one may also think of optimizations of
the algorithm itsself.
For example, the iteration
\bea\label{Eiterimp}
X_{n+1} = X_n - \frac{1}{2}\left[A e^{-X_n} - A^{-1} e^{+X_n}\right]
\eea
has in the 1-dimensional case near its fixed point $X_\infty=\ln A$ better
convergence properties than (\ref{Eiteration}).
However, the iteration does not converge far away from the solution.
Therefore, it could be used as an optional last refinement step after
the conventional iteration (\ref{Eiteration}) has converged with sufficient accuracy.

%%%%%%%%%%%%%%%%%%%%%%%%%%%%%%%%%%%%%%%%%%%%%%%%%%%%%%%%%%%%%%%%%%%%%%%%%%%%%%%
%%%%%%%%%%%%%%%%%%%%%%%%%%%%%%%%%%%%%%%%%%%%%%%%%%%%%%%%%%%%%%%%%%%%%%%%%%%%%%%
%%%%%%%%%%%%%%%%%%%%%%%%%%%%%%%%%%%%%%%%%%%%%%%%%%%%%%%%%%%%%%%%%%%%%%%%%%%%%%%
%%%%%%%%%%%%%%%%%%%%%%%%%%%%%%%%%%%%%%%%%%%%%%%%%%%%%%%%%%%%%%%%%%%%%%%%%%%%%%%

\section{A Numerical Algorithm}\label{Salgorithm}

The fixed-point iteration (\ref{Eiteration}) 
requires the calculation of the exponential of the iterates.
The associated computational burden can be reduced by exploiting that the proposed 
fixed-point iteration produces a series of mutually commuting matrices, if initialized properly.
Therefore, two successive exponentials can also be computed iteratively, which has 
the advantage that the norm of the matrix to be exponentiated in each step does not become
too large.
In this case, the inverse scaling and squaring method, which is based on
\bea\label{Eexponential}
e^B = \left[\exp\left(B/2^j\right)\right]^{2^j} 
\eea
is known to produce good results with a modest number of matrix multiplications 
\cite{higham2005a,moler2003a}.
For example using a $k$-th order Taylor approximant the exponential (\ref{Eexponential}) 
can be calculated with just $k+j-1$ matrix multiplications.
The algorithm can be summarized as follows
\begin{itemize}
\item Determine an initial matrix $X_0$
\begin{itemize} 
\item with eigenvalues close to those of $\ln A$
\item with $\commute{X_0}{A}=0$
\end{itemize}
\item set $Y_0 = \exp(-X_0)$
\item iterate
\bea
X_{n+1} &=& A Y_n - \f{1} + X_n\nn
Y_{n+1} &=& Y_n \exp\left\{-\left(A Y_n - \f{1}\right)\right\}\nonumber
\eea
until convergence is reached\\
(e.g. $\norm{X_{n+1}}\norm{X_{n+1}-X_n} \le \varepsilon$)\\
or a maximum number of iterations has been exceeded
\item refinement step [optional]:
\bea
X_{\rm fin} = X_n - \frac{1}{2}\left[A^{-1} e^{X_n} - A e^{-X_n}\right]\nonumber
\eea
\end{itemize}
Note that $Y_n$ will converge to the inverse of $A$ (although there exist
by far more efficient methods to achieve this).

%%%%%%%%%%%%%%%%%%%%%%%%%%%%%%%%%%%%%%%%%%%%%%%%%%%%%%%%%%%%%%%%%%%%%%%%%%%%%%%
%%%%%%%%%%%%%%%%%%%%%%%%%%%%%%%%%%%%%%%%%%%%%%%%%%%%%%%%%%%%%%%%%%%%%%%%%%%%%%%
%%%%%%%%%%%%%%%%%%%%%%%%%%%%%%%%%%%%%%%%%%%%%%%%%%%%%%%%%%%%%%%%%%%%%%%%%%%%%%%
%%%%%%%%%%%%%%%%%%%%%%%%%%%%%%%%%%%%%%%%%%%%%%%%%%%%%%%%%%%%%%%%%%%%%%%%%%%%%%%

\section{Performance Analysis}\label{Sperformance}

\begin{figure}
\includegraphics[height=7cm]{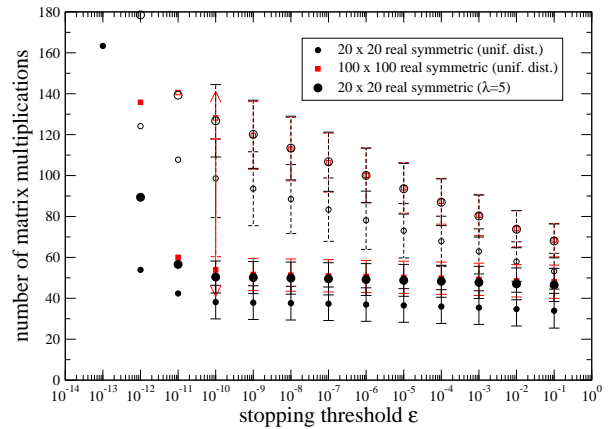}
\caption{\label{Fmatmul}Number of matrix multiplication required to achieve 
convergence to $\ln A$ for different stopping thresholds $\varepsilon$.
Hollow symbols include all matrix multiplications (including calculation Taylor approximants),
whereas full symbols refer to the multiplications required by the fixed point iteration scheme.
Error bars correspond to one standard deviation sampled over 1000 random matrices.
For the chosen test matrices, the algorithm shifts the computational burden towards
matrix exponentiation, and the total number of matrix multiplications
increases approximately linearly with the required precision.
Around $\varepsilon\approx 1\cdot 10^{-10}$ the scaling breaks down due to numerical roundoff
errors (error bars omitted), compare also figure \ref{Ferror}.
It is also visible that the dependence on the dimension of the matrix is rather weak.
}
\end{figure}

In order to estimate the performance of the algorithm, some sample matrices
have been generated.
For different matrix dimensions, 1000 matrices have been randomly generated.
Some test matrices had a uniform eigenvalue distribution in the interval $[1\cdot 10^{-8},1]$,
others were exponentially distributed according to 
$\rho_\lambda(x) = \lambda e^{-\lambda x}$ also in the interval $[1\cdot 10^{-8},1]$.
The diagonal matrix generated by these eigenvalues has been transformed into a 
non-diagonal test matrix by applying random orthogonal transformations $A=Q^{\rm T} A_{\rm D} Q$ with
$Q_{ij} = \delta_{ij} - 2 v_i v_j$ and $\sum_i v_i^2 = 1$.
All iterations were initialized with $X_0 = 2 A - [1+\ln(2)]\f{1}$, which is not
necessarily the optimum choice.
For all norm calculations, the Frobenius matrix norm has been used.
The iteration for the logarithm used as a stopping criterion 
$\norm{X_n}\norm{X_{n+1}-X_n} \le \varepsilon$, and the calculation of the 
exponential of a matrix $B$ used $\norm{B}^n \le \varepsilon$.
To calculate the matrix exponential, the scaling and squaring method was used in
combination with truncated Taylor approximants \cite{higham2005a}.
Note that the efficiency of the scaling and squaring method can in principle 
be increased by approximately 50\% if instead of Taylor approximants, Pad\'e 
approximations are used \cite{higham2005a}.
The number of matrix multiplications to obtain convergence was therefore counted 
with and without including those required by computing the Taylor approximants to the
matrix exponential, see figure \ref{Fmatmul}.
Whereas the number of total required matrix multiplications increases approximately linearly
with $\varepsilon$, it is visible that the algorithm shifts the computational burden towards
matrix exponentiation, since the number of remaining matrix 
multiplications is approximately constant.
The total number of matrix multiplications is competative with much more sophisticated 
algorithms existing in the literature \cite{cheng2001a}.

\begin{figure}
\includegraphics[height=7cm]{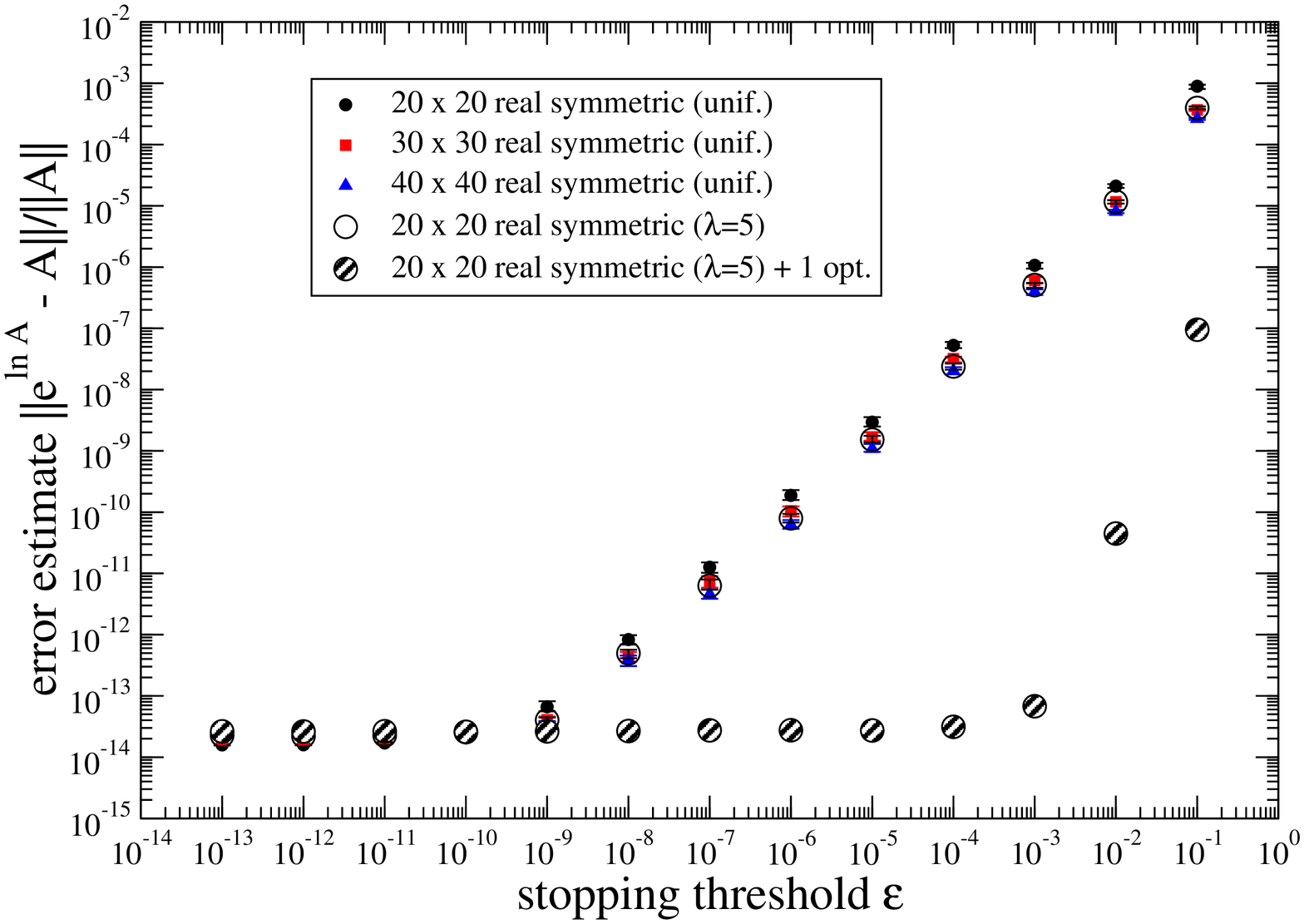}
\caption{\label{Ferror}
Scaling of the error estimate depending on the stopping criterion $\varepsilon$.
Each point displays the median (out of 1000 random samples) and error bars give the 
95\% confidence bounds for the median.
Above a threshold criterion of $\varepsilon \ge 1 \cdot 10^{-10}$, the error estimate 
scales approximately linearly with $\varepsilon$.
Applying only a single iteration of (\ref{Eiterimp}) with $\epsilon_{\rm final} = 10^{-16}$, 
the final error estimate can be reduced extremely (patterned symbols).
For small $\varepsilon$, a saturation due to finite machine precision is visible.
}
\end{figure}

In order to get an estimate on the error of $\ln A$ one can perform the inverse operation, i.e., 
exponentiate the result (with a much better precision) and compare it with the original
test matrix.
Thus, from $\norm{\exp(\ln A+\Delta) - A}/\norm{A}$ one has an estimate for the 
actual error $\norm{\Delta}$, see figure \ref{Ferror}.
It becomes visible that for $\varepsilon \ge 10^{-10}$ the accuracy of the computed 
logarithm is much better than one might have expected.
In addition, with a single application of (\ref{Eiterimp}), the final error of the
algorithm can be reduced by orders of magnitude.

%%%%%%%%%%%%%%%%%%%%%%%%%%%%%%%%%%%%%%%%%%%%%%%%%%%%%%%%%%%%%%%%%%%%%%%%%%%%%%%
%%%%%%%%%%%%%%%%%%%%%%%%%%%%%%%%%%%%%%%%%%%%%%%%%%%%%%%%%%%%%%%%%%%%%%%%%%%%%%%
%%%%%%%%%%%%%%%%%%%%%%%%%%%%%%%%%%%%%%%%%%%%%%%%%%%%%%%%%%%%%%%%%%%%%%%%%%%%%%%
%%%%%%%%%%%%%%%%%%%%%%%%%%%%%%%%%%%%%%%%%%%%%%%%%%%%%%%%%%%%%%%%%%%%%%%%%%%%%%%

\section{Summary}\label{Ssummary}

The fixed-point iteration presented here has some advantages in comparison to 
standard algorithms \cite{cheng2001a,cardoso2003a}.
One of the most important advantages is the ease of implementation.
Given the simplicity of the inverse-scaling and squaring algorithm for computing
the matrix exponential, the user basically has to implement matrix-matrix multiplication.
As has been demonstrated, the fixed point iteration guarantees fast convergence if
a good initial guess is given.
For some special test problems, the convergence is competative with state of the
art algorithms.
With some knowledge on the spectrum of the matrix of which the logarithm should
be taken, the convergence can be tuned to the specific problem.
In any case, the presented iteration can be used as a final step with a 
low-precision result produced with established methods \cite{cheng2001a,cardoso2003a}.

In the future, it would be interesting whether the iteration algorithm can be adapted
to increase convergence for a bad initial guess.

%%%%%%%%%%%%%%%%%%%%%%%%%%%%%%%%%%%%%%%%%%%%%%%%%%%%%%%%%%%%%%%%%%%%%%%%%%%%%%%
%%%%%%%%%%%%%%%%%%%%%%%%%%%%%%%%%%%%%%%%%%%%%%%%%%%%%%%%%%%%%%%%%%%%%%%%%%%%%%%
%%%%%%%%%%%%%%%%%%%%%%%%%%%%%%%%%%%%%%%%%%%%%%%%%%%%%%%%%%%%%%%%%%%%%%%%%%%%%%%
%%%%%%%%%%%%%%%%%%%%%%%%%%%%%%%%%%%%%%%%%%%%%%%%%%%%%%%%%%%%%%%%%%%%%%%%%%%%%%%
\section{Acknowledgements}

The author gratefully acknowledges financial support by the DFG 
grant \# SCHU~1557/1-2.

$^*${\tt schaller@theory.phy.tu-dresden.de}

%%%%%%%%%%%%%%%%%%%%%%%%%%%%%%%%%%%%%%%%%%%%%%%%%%%%%%%%%%%%%%%%%%%%%%%%%%%%%%%
%%%%%%%%%%%%%%%%%%%%%%%%%%%%%%%%%%%%%%%%%%%%%%%%%%%%%%%%%%%%%%%%%%%%%%%%%%%%%%%
%%%%%%%%%%%%%%%%%%%%%%%%%%%%%%%%%%%%%%%%%%%%%%%%%%%%%%%%%%%%%%%%%%%%%%%%%%%%%%%
%%%%%%%%%%%%%%%%%%%%%%%%%%%%%%%%%%%%%%%%%%%%%%%%%%%%%%%%%%%%%%%%%%%%%%%%%%%%%%%

%%%%%%%%%%%%%%%%%%%%%%%%%%%%%%%%%%%%%%%%%%%%%%%%%%%%%%%%%%%%%%%%%%%%%%%%%%%%%%%
%%%%%%%%%%%%%%%%%%%%%%%%%%%%%%%%%%%%%%%%%%%%%%%%%%%%%%%%%%%%%%%%%%%%%%%%%%%%%%%
%%%%%%%%%%%%%%%%%%%%%%%%%%%%%%%%%%%%%%%%%%%%%%%%%%%%%%%%%%%%%%%%%%%%%%%%%%%%%%%
%%%%%%%%%%%%%%%%%%%%%%%%%%%%%%%%%%%%%%%%%%%%%%%%%%%%%%%%%%%%%%%%%%%%%%%%%%%%%%%
%%%%%%%%%%%%%%%%%%%%%%%%%%%%%%%%%%%%%%%%%%%%%%%%%%%%%%%%%%%%%%%%%%%%%%%%%%%%%%%
%%%%%%%%%%%%%%%%%%%%%%%%%%%%%%%%%%%%%%%%%%%%%%%%%%%%%%%%%%%%%%%%%%%%%%%%%%%%%%%
%%%%%%%%%%%%%%%%%%%%%%%%%%%%%%%%%%%%%%%%%%%%%%%%%%%%%%%%%%%%%%%%%%%%%%%%%%%%%%%
%%%%%%%%%%%%%%%%%%%%%%%%%%%%%%%%%%%%%%%%%%%%%%%%%%%%%%%%%%%%%%%%%%%%%%%%%%%%%%%
%%%%%%%%%%%%%%%%%%%%%%%%%%%%%%%%%%%%%%%%%%%%%%%%%%%%%%%%%%%%%%%%%%%%%%%%%%%%%%%
%%%%%%%%%%%%%%%%%%%%%%%%%%%%%%%%%%%%%%%%%%%%%%%%%%%%%%%%%%%%%%%%%%%%%%%%%%%%%%%
%%%%%%%%%%%%%%%%%%%%%%%%%%%%%%%%%%%%%%%%%%%%%%%%%%%%%%%%%%%%%%%%%%%%%%%%%%%%%%%
%%%%%%%%%%%%%%%%%%%%%%%%%%%%%%%%%%%%%%%%%%%%%%%%%%%%%%%%%%%%%%%%%%%%%%%%%%%%%%%
%%%%%%%%%%%%%%%%%%%%%%%%%%%%%%%%%%%%%%%%%%%%%%%%%%%%%%%%%%%%%%%%%%%%%%%%%%%%%%%
%%%%%%%%%%%%%%%%%%%%%%%%%%%%%%%%%%%%%%%%%%%%%%%%%%%%%%%%%%%%%%%%%%%%%%%%%%%%%%%
%%%%%%%%%%%%%%%%%%%%%%%%%%%%%%%%%%%%%%%%%%%%%%%%%%%%%%%%%%%%%%%%%%%%%%%%%%%%%%%
%%%%%%%%%%%%%%%%%%%%%%%%%%%%%%%%%%%%%%%%%%%%%%%%%%%%%%%%%%%%%%%%%%%%%%%%%%%%%%%
\end{document}